\documentclass{ojac}
\usepackage{amsmath,amstext,amsthm,amssymb,amsxtra}
\usepackage[colorlinks,pagebackref,hypertexnames=false]{hyperref}
\usepackage[backrefs]{amsrefs}

\theoremstyle{plain} 
\newtheorem{lemma}[equation]{Lemma}
\newtheorem{proposition}[equation]{Proposition}
\newtheorem{theorem}[equation]{Theorem}

\theoremstyle{definition}

\theoremstyle{remark}
\newtheorem{remark}[equation]{Remark}
\newtheorem*{Acknowledgment}{Acknowledgment}

\numberwithin{equation}{section}

\def\norm#1.#2.{\lVert#1\rVert_{#2}}
\def\Norm#1.#2.{\bigl\lVert#1\bigr\rVert_{#2}}
\def\NOrm#1.#2.{\Bigl\lVert#1\Bigr\rVert_{#2}}
\def\NORm#1.#2.{\biggl\lVert#1\biggr\rVert_{#2}}
\def\NORM#1.#2.{\Biggl\lVert#1\Biggr\rVert_{#2}}

\def\ip#1,#2,{\langle #1,#2\rangle}
\def\Ip#1,#2,{\bigl\langle#1,#2\bigr\rangle}
\def\IP#1,#2,{\Bigl\langle#1,#2\Bigr\rangle}

\def\mid{\,:\,}

\def\abs#1{\lvert#1\rvert}
\def\Abs#1{\bigl\lvert#1\bigr\rvert}
\def\ABs#1{\biggl\lvert#1\biggr\rvert}

\def\XXint#1#2#3{{\setbox0=\hbox{$#1{#2#3}{\int}$}
     \vcenter{\hbox{$#2#3$}}\kern-.5\wd0}}

\def\eqdef{\stackrel{\mathrm{def}}{{}={}}}

\def\mid{\;|\;}

\begin{document}
\title{On an Argument of Shkredov on \\ Two-Dimensional Corners}

\author{Michael T. Lacey  
\\   
School of Mathematics\\
Georgia Institute of Technology\\
Atlanta GA 30332\\
\texttt{lacey@math.gatech.edu}
\and
William McClain
School of Mathematics\\
Georgia Institute of Technology\\
Atlanta GA 30332\\
\texttt{bill@math.gatech.edu}
}

\maketitle

\begin{abstract} Let $\mathbb F _2^n$ be the finite field of cardinality $2 ^{n}$.
For  all  large $n$, any subset $A\subset \mathbb F _2^n\times \mathbb F _2 ^n$
of cardinality
\begin{equation*}
\abs{ A} \gtrsim 4^n  \frac {\log \log n} {\log n}\,,
\end{equation*}
must contain three points $ \{(x,y)\,,(x+d,y)\,,(x,y+d)\}$ for
$x,y,d\in \mathbb F_2^n$ and $d\neq0$.   Our argument is an elaboration
 of an argument of Shkredov \cite {MR2223244},
building upon the finite field analog of Ben Green \cite {MR2187732}.
The interest in our result is in the exponent on $ \log n$, which is larger than
has been obtained previously.
\end{abstract}

\section{Main Theorem}
We are interested in  extensions of Roth's theorem on arithmetic
progressions in dense sets of integers to a two-dimensional,
finite-field setting. Specifically, for a finite group $ G$ define
the quantity $r_{\angle}(G)$ to be the cardinality of the largest
subset of $G\times G$ containing no \emph{corner}.  A \emph{corner} is
triple of points of the form $\{(x,y), (x + d,y), (x, y+d)\}$ with  $d
\neq 0$.

While this concept is most interesting in the context of the groups $G=\mathbb Z _{N}$,
it already makes sense--and is substantial--in the context of finite fields.
In this paper, we only consider the case of $ G= \mathbb F_2^n$.
Here and throughout this paper we write $N=2^n=\abs{\mathbb F_2^n}$.

\begin{theorem} \label{mainthm} $r_{\angle}(\mathbb{F}_2^n) \ll
N^2 \frac {\log \log \log N} {\log \log N}$.
\end{theorem}

This bound is an improvement, in the setting of $\mathbb F _2 ^{n}$,
of the bounds provided by Shkredov \cites{MR2223244,math.NT/0503639},
and as simplified by Ben Green \cites{MR2187732,green-onshk}.
The main point is that we elaborate on the
`Density Increment' procedure, obtaining a density increment on a set
which is the intersection of sublattices in two distinct sets of bases.

Our theorem is an example of the quantitative bounds on questions of arithmetic
combinatorics. We refer the reader to the papers of Gowers \cite {MR1844079}, and surveys by
T.~Tao \cite {math.NT/0505402} and Ben Green \cites{MR2187732,math.NT/0508063} for more history.

Erd\H os and Graham  raised the question of quantitative bounds for $r _{\angle}$, and this question
 was raised again by Gowers \cite {MR1844079}.
 Ajtai and Szemer{\'e}di \cite {MR0369299}
 first proved that $r _{\angle}(\mathbb Z _N)=o(N ^{2})$.
Furstenberg and Katznelson  \cites{MR1191743,MR833409} gave a far reaching
extension, though their method of proof does not in and of itself permit  explicit bounds.
Solymosi \cite {MR2047239},  and V.~Vu \cite {vu} provided such bounds, although of
a  weak nature.

I.~Shkredov \cites {MR2223244,MR2157918} provided the first `reasonable' bounds.  We are
using his ingenious argument, as explained and simplified by Ben Green \cite {MR2187732}
in the finite field setting. In particular Green showed that one
could achieve an estimate in which $ r _{\angle} (N)/N ^2 $
decreased like   $ (\log\log N) ^{-c}$ where $ c$ could be taken to be $ 1/21$.
We find an additional extension of this argument, and
sharpen some inequalities to obtain our Theorem.

We comment that `the finite field thesis' holds that questions of this type
should first be studied in the context of finite fields.  This is because one can implement
many of the tools of analysis, e.\thinspace g.~convolution and  Fourier transform, in that setting.
In addition, one has the powerful concept of being able to pass to appropriate affine
subspaces.  Moreover, there are
a range of methods that one can use to `lift' the finite field argument to $\mathbb Z _{N}$.
See papers by Bourgain \cite {MR1726234}, Green and Tao \cite {math.NT/0503014} and
Shkredov \cite {math.NT/0503639} 
for more information.

\begin {Acknowledgment} 
One of us (M.L.) was supported in part by the NSF and the Guggenheim Foundation. 
One of us (W.M.) was supported in part by a VIGRE grant awarded to Georgia 
Institute of Technology.  Both authors are grateful for the hospitality of
the University
of British Columbia. While there, we benefited greatly from conversations with
Izabella {\L}aba and Jozef Solymosi. Ben Green and I.~Shkredov pointed out an error in
an earlier version of this paper.  We also thank the referee for additional 
corrections, which lead to strengthening of the main result. 
\end{Acknowledgment}

\subsection{Outline of Proof}

Ben Green has provided a comprehensive
outline of the method of proof \cite {MR2187732}, so
we will be somewhat brief.

Let $A\subset \mathbb F _2^n\times \mathbb F _2 ^{n}$.
It is natural to count the
expected number of corners that $A$ has.
This expectation  will be approximately
what it should be if the `Box norms' of $A$ are small.
These norms, we use three of them,
 are  explicitly given in (\ref{e.box})---(\ref{e.BoxNormYD}), and
are a two dimensional analog of the `combinatorial square' norms
that play such a prominent role in the proof of Roth's Theorem
\cites {MR1335785
,
MR2187732
,
math.NT/0503014}.  It is therefore a certain measure of `uniformity.'

An important difference between our approach and the previous ones is that we emphasize the
role of three distinct coordinates in the problem. Two of these are the obvious
$ X$ and $ Y$ coordinates, given by the canonical basis $ \operatorname e_1$
and $ \operatorname e_2$ respectively.  The third is the \emph{diagonal} coordinate $ D$,
given by $ \operatorname e_1+\operatorname e_2$.
In this way, our argument resembles the ergodic theoretic 
arguments, and in particular that of Conze and Lesigne \cite{MR788966}
who describe the characteristic factor for the ergodic averages\footnote{
The argument for two commuting transformations in this paper is complete. We are
grateful to Bryna Kra for pointing out this reference to us and to Emmaunel Lesigne
for providing us with a copy of this paper.}
\begin{equation*}
N ^{-1} \sum _{n=1} ^{N} f_1 (T_1 ^{n}) f _2(T_2 ^{n})
\end{equation*}
where $ T_1,T_2$ are commuting measure preserving transformations.

Two of the three  box  norms employ this diagonal coordinate.
Either being large is an `obstruction to uniformity.' (See \cite {math.NT/0505402}.)
But, in contrast to the case of Roth's theorem, this obstruction to uniformity has
no clear arithmetic consequence.
It does imply, however,  that
the set $A$ has an increased density on a sublattice.  That is, there is a subsets
$X, Y, D\subset \mathbb F _2^n$,  of relatively large density,
for which $A$ has a larger density on the \emph{intersection of
$ X\times Y$ and $ X\stackrel{\textup{diag}}{\times} D$, }
where in the latter product, we are  taking the product in the coordinate system
$ (\operatorname e_1, \operatorname e_1+\operatorname e_2)$.
See Lemma~\ref{increment}.    Forming the intersection of two sublattices in this way
is the main new contribution of this paper.

It was a significant insight of Shkredov  \cite {MR2223244} that (1) one could, after an additional
argument, assume that $X$, $ Y$ and $ D$
had an arithmetic structure, namely that it was \emph{ uniform}
in the sense of (\ref{e.Uni}), see Lemma~\ref{uniformity};  and (2) the Box norm could still be used
as an `obstruction to uniformity' if $X$, $ Y$  and $ D$ were uniform.  This is the content of the
Lemma~\ref{gvn}.

If there is an obstruction to uniformity, an increment in the density of $A$ can be found.
One can find the obstruction to uniformity only a finite number of times, else the
density of $A$ would exceed one.  Thus at some point, there would be no obstruction
to uniformity, and so $A$ would have a corner.  All of the details of the proof
are below.

\section{Preliminaries and Definitions}

We will let $H\subset \mathbb F _2^n$ denote a subspace.
The dimension of this subspace will decrease in the course of the proof.
$X, D, Y\subset H$ denote  subsets.  We adopt the notations of
probability and expectation with respect to the counting measure
on $H$.  We view $ X$ as a subset of the first coordinate associated
to basis element $ \operatorname e_1$; $ Y$ as a subset
of the second coordinate, associated to basis element $ \operatorname e_2$;
and $ D$ as a subset of the `diagonal' coordinate associated to $ \operatorname e_1+
\operatorname e_2$.
We will be working with subsets of
\begin{equation}\label{e.Sdef}
S \eqdef X\times Y\cap X\stackrel{\textup{diag}}{\times} D.
\end{equation}
It is to be emphasized that the products are taken in the first instance in the
$ \operatorname e_1$ and $ \operatorname e_2$ coordinates; the second in the $ \operatorname e_1$
and $\operatorname e_1+\operatorname e_2$ coordinates.

Write \emph{the density of $ X$
in $ H$} as
\begin{equation*}
\delta_X := \mathbb P (X\mid H)= \tfrac {\abs{ X}} {\abs{H}}\,,
\end{equation*}
and similarly for  $ \delta _{Y}$ and $ \delta _D$.
In the iterative procedure used in this proof, these densities will
be decreasing, somewhat rapidly.
Throughout this paper we will not only be concerned with the density
of the subsets $X$ and $ D$, but also with how `uniformly distributed' they are
in $H$.  A quantification of this quality comes in the following
definition:
\begin{equation} \label{e.Uni}
\norm X.\textup{Uni}. := \sup_{\xi \neq 0} \tfrac
{\abs{\widehat{X}(\xi)}} {\abs{H}}\,.
\end{equation}
This definition only makes sense relative to a subspace $ H$.
If $\norm X.\textup{Uni}. \le \eta$ then we say that $X$ is
\emph{$\eta$--uniform}.
Here $\widehat{X}$ represents the the Fourier transform of $X$ which is
defined as follows, in which $\omega= {-1}$ is a second root of unity:
\begin{equation*}
\widehat{g}(\xi) := \sum_{x \in H} g(x) \omega^{x \cdot \xi}.
\end{equation*}

It is immediate from the definition that a translate of a uniform set by
an element of $ H$  is again uniform.
After the deletion of a small subset, a uniform set is again uniform.  Let
$ E\subset X$, we have
\begin{equation}\label{e.thin}
\norm X-E.\textup{Uni}.\le \norm X.\textup{Uni}.+\mathbb P (E\mid X)\,.
\end{equation}

This proposition will be used repeatedly.

\begin{proposition}\label{uniformconvolution}
Let $X \subset H$ and denote the density of $X$ in $H$ by
$\delta_X$.  Then
\begin{equation}\label{uniconv}
\bigl[\mathbb E_{d \in H} \Abs{ \mathbb E_{y \in H} X(d-y)G(y) -
\delta_X  \mathbb E_{y \in H} G(y) } ^2\bigr] ^{1/2}  \le \norm X. \textup{Uni}.
\bigl[ \mathbb E _{y}G (y )^2  \bigr] ^{1/2}
\end{equation}
for any function G.
\end{proposition}

Notice that we are  comparing a  convolution  to it's zero Fourier mode. The Proposition
follows from Plancherel, and the definition of uniformity.  A form of this inequality
that we will use several times is this: For any two sets $ A,B\subset H$,
\begin{equation}\label{e.justone}
\mathbb E _{x\in X}\Abs{ \mathbb E _{\in Y} A(x+y)B (y)- \mathbb P (A) \mathbb P (B)} ^2
\le\min \{ \norm A.\textup{Uni}. ^{1/2}  \mathbb P (B)\,,\,
\norm B.\textup{Uni}. ^{1/2}  \mathbb P (A)\}\,.
\end{equation}
That is, only uniformity in one coordinate is required.


Now, $ S$ is as in (\ref{e.Sdef}), and  let $A \subset S$. Write \emph{the density of A} as
\begin{equation*}
\delta :=\mathbb P(A\mid S)\,.
\end{equation*}
This quantity will increase in the iterative procedure used in the proof.
We define the \emph{balanced function} of $A$ to be the function
supported on  $S$ as
\begin{equation*}
f(x,y) := A(x,y) - \delta S\,.
\end{equation*}

Our standing assumption will be
\begin{equation}\label{e.upsilon}
\begin{split}
\norm X.\textup{Uni}.\,,\, \norm Y.\textup{Uni}.\,,\,
\norm D.\textup{Uni}.\le \upsilon \,,
\qquad
\upsilon \eqdef (\delta \delta _X \delta _Y \delta _D) ^{C}
\end{split}
\end{equation}
where $ C $ is a large constant which we need not specify exactly, as its precise
value only influences implied constants in our main Theorem.
In the proofs of Lemmas, we will use the notation $ \upsilon '$ for
a fixed, but unimportant, function of $ \upsilon $, that tends to zero as $ \upsilon $ does.

Further, for a function $ f\,:\, S \longrightarrow \mathbb C
$, define the following norm
\begin{equation}\label{e.box}
\norm f.\Box. \eqdef \delta ^{-4} _{D}
\mathbb E _{\substack{x,x'\in X\\ y,y'\in Y }}
f (x,y)f( x',y) f (x,y') f (x',y')
\end{equation}
where we use the standard basis $ (\operatorname e_1,\operatorname e_2)$.
This  norm averages `cross correlations' of $ f$ over all boxes
in $ X \times Y $.  When $ f$ is the balanced function
of $ A$, the norm being `large' is an obstacle to $ A$ having the correct number
of corners.

We use two additional norms. In the
$ (\operatorname e_1,
\operatorname e_1+\operatorname e_2)$ coordinate system,
\begin{equation}\label{e.BoxNormXD}
\norm f. \Box,X. ^{4}:=   \delta _{Y} ^{-4}\mathbb E _{\substack{ x, x'\in X\\
 d, d'\in D}} f( x, d)f( x', d)f( x, d')
f( x', d')\,.
\end{equation}
Similarly, with respect to the $ (\operatorname e_2,
\operatorname e_1+\operatorname e_2)$ coordinate system, define
\begin{equation}\label{e.BoxNormYD}
\norm f. \Box,Y. ^{4}:=   \delta _{X} ^{-4}\mathbb E _{\substack{ y, y'\in Y\\
 d, d'\in D}} f( y, d)f( y', d)f( y, d')
f( y', d')\,.
\end{equation}

The leading normalizations in these definitions are initially confusing, but chosen
so that the norms are essentially bounded from above by the $ L^\infty $ norm of $ f$.
The point of these next propositions is that the quantities introduced above behave
as they should, under the assumption of uniformity. In particular,
(\ref{e.BoxCount1}) and (\ref{e.BoxCount2}) justify the normalizations
in the definition of the Box norms.

\begin{proposition}\label{p.upsilon1}
Let $ X,Y,D\subset H$ be as above, and let $   S$
be as in (\ref{e.Sdef}).
Assuming (\ref{e.upsilon}) we have
\begin{gather}\label{e.Scount1}
\mathbb E _{\substack{x\in X\\ y\in Y }}S=\delta _{D} +O(\upsilon') ;
\\
\label{e.Scount2}
\mathbb E _{\substack{x\in X\\ d\in D }}S=\delta _{Y}+O(\upsilon')  ;
\\
\label{e.Scount3}
\mathbb E _{\substack{d\in D\\ y\in Y }}S=\delta _{X}+O(\upsilon'   );
\\
\label{e.Tcount}
\mathbb E _{\substack{x \in X\\ y\in Y \\ s\in H}}
S (x,y)
S (x+s,y) S (x,y+s)
=\delta _X  \delta _Y   \delta _D^2 +O(\upsilon');
\\
\label{e.BoxCount1}
\mathbb E _{\substack{ x, x'\in X\\
 d, d'\in D}} S( x, d)S( x', d)S( x, d')
S( x', d')=\delta _{Y} ^{4}+O(\upsilon');
\\  \label{e.BoxCount2}
\mathbb E _{\substack{ y, y'\in Y\\
 d, d'\in D}} S( y, d)S( y', d)S( y, d')
S( y', d')=\delta _{X} ^{4}+O(\upsilon').
\end{gather}
In (\ref{e.BoxCount1}) we are using the
$ (\operatorname e_1, \operatorname e_1+\operatorname e_2)$ coordinate systems,
while in the (\ref{e.BoxCount2}) we are using the
$ (\operatorname e_2, \operatorname e_1+\operatorname e_2)$
 coordinate system.
\end{proposition}

\begin{proof}
For (\ref{e.Scount1}), observe that
\begin{align*}
\delta _X \delta _Y \mathbb E _{\substack{x\in X\\ y\in Y }}S
 &= \mathbb E _{x,y\in H} X(x)Y(y)D(x+y)=\delta _X \delta _Y \delta _D+O( \delta _X ^{1/2}
 \delta_Y ^{1/2}\upsilon) \,.
\end{align*}
This equality only  requires uniformity in one coordinate.
See (\ref{e.justone}).  The equalities
(\ref{e.Scount2}) and (\ref{e.Scount3}) are  corollaries, after a change of variables.

\smallskip

To see (\ref{e.Tcount}), we apply Lemma~\ref{l.BeyondConvolution}, with
$ f=X$ and $ g=Y$.  Using the notation $ \Phi   $ in (\ref{e.phi}), we have
\begin{align*}
\mathbb E _{x,y,s\in H} & S (x,y) S (x+s,y) S (x,y+s)
\\
&= \mathbb E _{x,y,s\in H} X (x)  X (x+s) D (x+y) D (s) Y (y) Y (y+s)
\\
&=
\mathbb E _{x,y\in H} X (x)  Y (y) D (x+y) \Phi (x+y)+O( \upsilon')
\\
&= \delta _X \delta _Y \mathbb E _{x,y\in H} D (x+y) \Phi (x+y)+O( \upsilon')
\\
&= \delta _X \delta _Y \mathbb E _{x\in H} D (x) \Phi (x)+O( \upsilon')
\end{align*}


It remains to estimate the last expectation, which we view as an inner product.
Observe that $ \widehat \Phi (0)= \abs{ H}  \delta _D \delta _X \delta _Y $. And,
by Plancherel,
\begin{align*}
\Abs{\mathbb E _{x\in H} D (x) \Phi (x) - \delta _D ^2 \delta _X \delta _Y}& =
	\abs{ H} ^{-2} \ABs{\sum _{\alpha \neq 0} \widehat D (\alpha )
	 \widehat \Phi (\alpha ) }
\\& \le
\norm D .\textup{Uni}.  \abs{ H}^{-1} \sum _{\alpha\neq0 }
\abs{  \widehat \Phi (\alpha ) }
\\ & \le
\upsilon  \abs{ H}^{-2} \delta _D  \sum _{\alpha\neq0 }
\abs{ \widehat X (\alpha )  \widehat Y (\alpha )}
\\ &\le
\upsilon  \delta _D \delta _X ^{1/2} \delta _Y ^{1/2} \,.
\end{align*}

\smallskip

Concerning (\ref{e.BoxCount1}), we use a similar proof to the one above.
We follow the notation in Lemma~\ref{l.BeyondConvolution}, and its proof.
Set
\begin{align*}
\Psi  (x,x') & \eqdef \mathbb E _{d} Y (x+d) D (d) Y (x'+d)\,,
\\
\Phi (x) & \eqdef \frac {\delta _D} {\abs{ H} ^2 }
\sum _{\alpha \in H} \widehat Y (\alpha ) ^2 \omega ^{\alpha \cdot x}\,.
\end{align*}
Lemma~\ref{l.BeyondConvolution} implies that $ \Psi $ is well approximated
by $ \Phi $.  Hence, we can estimate
\begin{align*}
\delta _{X} ^2  \delta _{D} ^2 \cdot
(\ref{e.BoxCount1}) & =  \mathbb E _{\substack{ x, x'\in H\\
 d, d'\in H}} X (x)  X (x') D (d) D (d')
 Y( x+ d)Y( x'+ d)Y( x+ d')
Y( x'+ d')
\\
&=
\mathbb E _{x,x'\in H} X (x) X (x' ) \Psi (x,x') ^2
\\
& = \mathbb E _{x,x'\in H} X (x) X (x' ) \Phi (x+x')^2 +O (\upsilon ')
\\
& =  \delta _X ^2 \mathbb E _{x\in H} \Phi (x+x') ^2 + O (\upsilon ')
\\
& = \delta _X ^2  \delta _D ^2  \abs{ H} ^{-4} \sum _{\alpha \in H}
\widehat Y (\alpha ) ^{4} + O( \upsilon ')
\\
&= \delta _X ^2 \delta _D ^2 \delta _Y ^{4} + O (\upsilon ')\,.
\end{align*}
Here, we have used Lemma~\ref{l.BeyondConvolution}, uniformity in $ X$,
Plancherel, and uniformity in $ Y$.  The second equality  (\ref{e.BoxCount2})
follows from the first. This completes the proof.

%

\end{proof}

\begin{remark}\label{r.Tcount} There is a second way to see (\ref{e.Tcount}), which
we only briefly indicate, since the method of proof is not self contained.
Consider $ S=X \times Y\cap X \stackrel {\textup{diag}} {\times } D$ as
a subset of $ X \times Y$, and let $ \Delta $ be it's balanced function.  Namely
$ \Delta =S - \mathbb P (S\mid X \times Y) X \times Y$.  One can then define
the Box norm, as does Shkredov
\begin{equation*}
\norm \Delta . \textup{RectBox}.  ^{4} \eqdef \mathbb E _{\substack{x,x'\in X\\ y,y'\in Y }}
\Delta (x,y) \Delta (x',y) \Delta (x,y') \Delta (x',y')
\end{equation*}
It follows from the proof of (\ref{e.BoxCount1}), that we have
$ \norm \Delta .\textup{RectBox}. \lesssim \upsilon '$.  Shkredov
\cite{MR2223244}
showed that under this assumption,
and uniformity in $ X$ and $ Y$, that the set $ S$ has nearly the expected number
of point in it.  That is the content of his `Generalized von Neumann Lemma.'
\end{remark}


\begin{lemma}\label{l.BeyondConvolution}  Let $ D $ be uniform, and
let $ f,g$ be two functions on $ H$.  Define
\begin{equation} \label{e.phi}
\Phi (x) \eqdef \frac{\delta _D} {\abs{ H} ^2 } \sum _{\alpha \in H}
\widehat f (\alpha ) \widehat g (\alpha ) \omega ^{\alpha \cdot x}\,.
\end{equation}
We have the inequality
\begin{equation}\label{e.BeyondConvolution}
\begin{split}
\Bigl[ \mathbb E _{x,y\in H} \Abs{\mathbb E _{s} &f (x+s) D (s) g (y+s) - \Phi (x+y)} ^2
\Bigr]^{1/2} 
\\&
\lesssim
\norm D.\textup{Uni}. \bigl[ \mathbb E _{x} f (x) ^2  \bigr] ^{1/2}  \cdot
\bigl[ \mathbb E _{y} g (y) ^2  \bigr] ^{1/2}
\end{split}
\end{equation}
\end{lemma}

\begin{proof}
Consider
\begin{equation*}
\Psi (x,y)=
\mathbb E _{s} f (x+s) D (s) g (y+s)
\end{equation*}
as a function on $ H \times H$.  Expanding $ f$ in dual variable $ \alpha $
and $ g$ in dual variable $ \beta $ we have
\begin{align*}
\Psi (x,y) & = \abs{ H} ^{-2} \sum _{\alpha ,\beta \in H}
 \widehat f (\alpha ) \widehat g (\beta ) \omega ^{
 \alpha \cdot  x+\beta \cdot y}
 \mathbb E _{s} D (s) \omega ^{(\alpha +\beta ) \cdot s}
 \\
 &= \abs{ H} ^{-2} \sum _{\alpha ,\beta \in H}
 \widehat f (\alpha ) \widehat g (\beta ) \omega ^{
 \alpha \cdot  x+\beta \cdot y}
 \frac {\widehat D (\alpha +\beta ) } {\abs{ H}}.
\end{align*}
This shows that $ \widehat \Psi (\alpha ,\beta ) =  \widehat f (\alpha ) \widehat g (\beta )
\abs{ H} ^{-1} \widehat D (\alpha +\beta ) $.
Clearly, $ \Phi  $ of the Lemma consists of
the reconstruction of $ \Psi $ from those Fourier coefficients
$ (\alpha ,\beta )$ for which $\alpha +\beta =0 $.
And by Plancherel, the Lemma follows from
\begin{align*}
 \sum _{\alpha +\beta \neq 0}  \abs{ \widehat \Psi (\alpha ,\beta ) } ^2
& \le \norm D. \textup{Uni}. ^2 \sum _{\alpha ,\beta } \abs{ \widehat f (\alpha ) \widehat g (\beta )}
^2
\\
&=\abs{ H} ^{4} \cdot \norm D. \textup{Uni}. ^2   \cdot
\mathbb E _{x} f (x) ^2     \cdot
 \mathbb E _{y} g (y) ^2  \,.
\end{align*}

\end{proof}

\section{Primary Lemmata}

Our first lemma is a generalized von Neumann estimate, a term coined
by Green and Tao
\cites{MR2187732,math.NT/0404188,math.NT/0503014}.  It gives
us sufficient conditions from which to conclude that A has a corner.

\begin{lemma}[Generalized von Neumann]\label{gvn}
Suppose that $ A\subset S$ with $ \mathbb P (A\mid S)=\delta $;
(\ref{e.upsilon}) holds; and  we have the two inequalities
\begin{gather} \label{trivialcorners}
 \delta _{X}  \delta _Y\delta _D \delta ^{2} N > C \,,
\\ \label{gvnbound}
\max\{ \norm f.\Box.\,,\,
\norm f.\Box,X.\,,\, \norm f.\Box,Y. \}\le \kappa \delta ^{5/4}\,.
\end{gather}
Then A has a corner.
\end{lemma}

Here, and throughout this paper, $C$ represents a large absolute
constant.  The exact value of $ C$ does not  impact the
qualitative nature of our estimate, so we do not seek to specify an
optimal value for it. Also $ 0<c,\kappa,\kappa '<1 $ are  small fixed constants,
which plays a role similar to $ C$.

Our second lemma tells us that if the conditions of our first lemma
are not satisfied then we can find a sublattice on which A has
increased density.

\begin{lemma}[Density Increment]\label{increment}
For $ 0<\kappa $ there is a constant  $ 0<\kappa ' <1$ for which the following holds.
Suppose that $ A\subset S= X \times Y\cap
X\stackrel{\textup{diag}}{\times} D$ with $ \mathbb P (A\mid S)=\delta $,
that $ f$ is the balanced function of $ A$ on $S$, and that
 \begin{equation*}
\max\{ \norm f.\Box.\,,\, \norm f.\Box,X. \,,\, \norm f.\Box,Y.\} >
\kappa\delta^{5/4} \,.
\end{equation*}
Then there exists $X' \subset X$, $ Y'\subset Y$, $ D'\subset D$
such that three conditions hold.
\begin{gather} \label{e.1outof3}
\textup{either $ X'=X$, or $ Y'=Y$, or $ D'=D$; }
\\
\label{denseinc}\mathbb P(A\mid S')
\ge \delta +\kappa '\delta ^{2 }  \,,\qquad
S'=X'\times Y'\cap X'\stackrel{\textup{diag}}{\times} D'\,;
\\
\label{newsize} \mathbb P (X'\mid X)
\,,\,\mathbb P (Y'\mid Y)
\,,\, \mathbb P (D'\mid D)\ge
\kappa'\delta^{2 }.
\end{gather}
\end{lemma}

We need only refine two of the three sets
$ X$, $ Y$ and $ D$ above. Note that with  uniformity in coordinate that is
\emph {not} refined, we
then have that the set $ S'=X'\times Y'\cap X'\stackrel{\textup{diag}}{\times} D'$
has about the expected number of points in it.

Our third lemma is a modification of one
from a note of Ben Green \cite {MR2187732}.  It
tells us that we can find a \emph{uniform} sublattice on which A has
increased density which is important since it is a required premise
in applying the generalized von Neumann Lemma.

\begin{lemma}[Uniformizing a Sublattice]\label{uniformity}
Suppose that   $ X,Y,D$ are as above.
In addition
\begin{enumerate}
\item $ X$, $ Y$, and $ D$
satisfy (\ref{e.upsilon});

\item  $ X'\subset X$, $ Y'\subset Y$, $ D'\subset D$, with $ \mathbb P (X'\mid
X)\ge c \delta ^{2} $ and similarly for $ Y$ and $ D$;

\item Either  $ X'=X$, $ Y'=Y$ or $ D'=D$;

\item  $ S'=X'\times Y'\cap X'\stackrel{\textup{diag}}{\times} D'$;

\item   $ \mathbb P (A\mid S')=\delta + c\delta ^{2 } $;

\item
$\operatorname {dim}(H)>C [ \delta ^{4 } (\upsilon'')^2 ] ^{-1} $, where $ 0<\upsilon''<1$ is fixed.

\end{enumerate}
Then there exists $X''\subset
X'$, $ Y''\subset Y$, $ D''\subset D'$ and $H', H''$, translates of the
\emph{same} subspace $ H_0\le H$, so that
\begin{gather}
\label{conc1}\norm X''.\textup{Uni}. \,,\, \norm Y''.\textup{Uni}.
\,,\, \norm D''. \textup{Uni}. \le \upsilon '',
\\
\label{conc2}\mathbb P(A\mid S'') \ge \delta + \tfrac
{c} {2}\, \delta^{2 }\,,\qquad S''=X''\times Y''\cap X''\stackrel{\textup{diag}}{\times} D''\,,
\\
\label{conc3}\operatorname {dim}(H_0) \ge \operatorname {dim}(H) -
C [ \delta ^{4} (\upsilon'')^2 ] ^{-1}\, ,
\\
\label{conc4}\mathbb P (X''\mid H')\ge \kappa  \delta ^{2 } \mathbb P (X'\mid H)\,.
\end{gather}
An inequality similar to the last one also holds for $ Y''$ and $ D''$.
In particular $ \mathbb P (D''\mid H'+H'')\ge \kappa  \delta ^{2 } \mathbb P (D'\mid H)$
\end{lemma}

 It is to be emphasized that
$ H'$ and $ H''$ are translates of the same subspace of $ H_0<H$,
therefore $ H'+H''$ is also a translate of $ H_0$.  Thus, after a joint
translation of $ A$, $ X$, $ Y$ and $ D$, we can assume that  $ H'$ and $ H''$ are
in fact the same subspace $ H$.  It is this translation that
is used in the iteration of the proof.

\section{Proof of Theorem \ref{mainthm}}

Combining Lemmas~\ref{gvn} through \ref{uniformity} of the previous
section yields the proof of Theorem~\ref{mainthm}.  Since the proof
is by recursion, we describe the conditional loop needed for the proof.

\begin{proof}

Initialize $X \leftarrow  \mathbb F_2^n$, $ Y\leftarrow \mathbb F _2 ^{n}
$, $ D\leftarrow  \mathbb F _2 ^{n}$,  $ S\leftarrow \mathbb F _2 ^{n}
\times \mathbb F _2 ^{n}$,  $H \leftarrow  \mathbb F_2^n$.
Likewise $\delta _X, \delta _Y, \delta _D\leftarrow 1$.
Fix a set $A_0$
with density $\delta _0$ in $\mathbb F ^n_2\times \mathbb F ^n_2$.  Initialize
$A\leftarrow A_0$ and $ \delta \leftarrow \mathbb P (A\mid S)$.  Now we
will iteratively apply the following steps:
\begin{enumerate}
\item  If $ \max\{ \norm f.\Box.\,,\,\norm f.\Box,X.\,,\, \norm f.\Box,Y.\}
> \kappa\delta^{5/4}$, apply Lemma~\ref{increment}.

\item If $X'$, $ Y'$ or $ D'$ is not $\upsilon =(\delta \delta_{X'} \delta _{Y'}
\delta _D) ^C $ uniform, apply Lemma~\ref{uniformity}. Suppose these sets are as in the Lemma:
subsets $ X''\subset X'$, $ Y''\subset Y'$, $ D''\subset D'$
and affine subspaces $H',H''\subset H $ containing $ X'',Y'',D''$.
After joint translation of $ X'', Y'', D'', A$ and
$ H', H''$, we can assume that  $ H'=H''$  and are  subspaces of $ H$.

\item Update variables:
\begin{gather*}
X \leftarrow X''\,,\quad Y\leftarrow Y''\,,
\quad D\leftarrow D'', \quad H \leftarrow H'
\,,
\\
\delta _X \leftarrow \mathbb P (X''\mid H')\,,
\quad \delta _Y \leftarrow \mathbb P (Y''\mid H')\,,
\quad \delta _D \leftarrow \mathbb P (D''\mid H')\,,
\\  S\leftarrow X \times Y \cap X \stackrel {\textup{diag}} \times  D\, ,
\quad \delta \leftarrow \mathbb P (A\mid S).
\end{gather*}

\item
Observe that the density of the incremented $A$ on the set $ S$
has increased  by at  least $\kappa \delta_0 ^{2 }   $. Also, the
incremented densities, $\delta _X$, $ \delta _Y$, and $ \delta _D$ have decreased by no more than
$(\kappa \delta_0)  ^{C}$.

\end {enumerate}

Once this loop stops, Lemma \ref{gvn} applies---provided that the initial dimension
is large enough---and we
conclude that A has a corner. This iteration  must stop in $ \lesssim \delta ^{-1 } _{0}$
iterates, else the density of $A$ on the sublattice would exceed one.
Thus we need to be able to apply Lemmas~\ref{increment} and \ref{uniformity}
$ \lesssim \delta_0 ^{-1 } $ times.
In order to do that, both $X$ and $H$ must be sufficiently large
at each stage of the loop.

This requirement places several lower bounds on $N=2^n$.
The most stringent of these comes from the loss of dimensions in
(\ref{conc3}).   Note that before the loop terminates, we can have
$\delta _X$  as small as
\begin{equation*}
\delta _X \ge (\kappa \delta_0 )^{ (\kappa \delta _0) ^{-1 } }\,.
\end{equation*}
 In order to apply Lemma~\ref{uniformity} at that stage, we need
\begin{equation*}
N >   2^{(C\delta_0)^{-C\delta_0 ^{-1 } }}.
\end{equation*}
From this condition we get the bound stated in the Theorem.
\end{proof}

\section{Proof of Generalized von Neumann Lemma}

\begin{proof}[Proof of Lemma~\ref{gvn}]

Define \begin{align}  \label{e.Tdef}
\operatorname T(f,g,h) & = \mathbb E_{x,s,y
\in H} f(x,y)g(x+s,y)h(x,y+s)\,.
\end{align}
Thus, $ \operatorname T (A,A,A)$ is the expected number of corners in $ A$.
 We show that this quantity
 is at least a fixed small multiple of $ \delta _{X} ^2 \delta _Y ^2
 \delta _D ^2 \delta ^{3}$.
 By assumption (\ref{trivialcorners}), $ C \delta _{X} ^2 \delta _Y ^2  \delta _D ^2
 \delta ^{3} N ^{3}> \delta \delta _X \delta _Y \delta _D N ^2$.
 The left hand side is the expected number of corners in $ A$, while the right is the
 number of trivial corners in $ A$---that is the number of points in $ A$.  Thus
 $ A$ is seen to have a corner.

 Throughout the proof, it is convenient to make the substitution
 $ s\to x+y+s$ in the expression for $ \operatorname  T (f,g,h)$, thus
 \begin{equation*}
\operatorname T (f,g,h)=\mathbb E _{x,y,s} f (x,y) g (y+s,y) h (x,x+s)\,.
\end{equation*}
 We are of course using the fact that we work in a field of characteristic two, but
 there is a similar substitution for any field.

 We make the substitution $A = f + \delta S$    to get
\begin{align} \label {0f}
\operatorname  T(A,A,A)&=  \delta ^{3}\operatorname T(S,S,S)
\\ \label{1f}
&\quad + \delta^2\operatorname T(f,S,S) +  \delta^2\operatorname T(S,f,S)
+ \delta^2\operatorname T(S,S,f)
\\ \label{2f}
&\quad + \delta\operatorname T(f,f,S) +  \delta\operatorname T(S,f,f)
+ \delta\operatorname T(f,S,f)
\\ \label{3f}
&\quad + \operatorname T (f,f,f)\,.
\end{align}
We have grouped the terms according to the number of $ f$'s that appear.

The main term is the right hand side of  (\ref{0f}).  Using (\ref{e.Tcount}), we see that
 $\delta ^{3}\operatorname T(S,S,S)\ge \tfrac 12
\delta _{X} ^2 \delta _Y ^2 \delta _D ^2 \delta ^{3}$.  (Observe the difference in the normalizations
on the expectations in (\ref{e.Tcount}).)


All three terms in (\ref{1f})
are approximately zero, but we have to use uniformity to see this.
For instance,   we appeal to (\ref{uniconv})
and (\ref{e.upsilon}) to see that
\begin{equation}\label{e.1f}
\begin{split}
\delta ^{2}\operatorname T (S,f,S) &=\delta ^2
\mathbb E _{ x, y, s\in H} X (x) Y ( x+ s) D (x+y)  f(y+s, y)
\\
&=\delta ^2 \delta _X
 \mathbb E _{ x, y, s\in H} Y ( x+ s)  D (x+y)  f ( y+ s,y) +O (\upsilon')
 \\
&=\delta ^2 \delta _X  \delta _Y \delta _D
 \mathbb E _{  y, s\in H}   f ( y,  s) +O (\upsilon')
 \\
&=O(\upsilon')\,.
\end{split}
\end{equation}
In this line and below, $ \upsilon '$ is an unimportant function of $ \upsilon $
which tends to zero.

The  three terms in (\ref{2f}) are all controlled by appeal to Lemma~\ref{l.T}.
For instance, we have using the assumption about the maximal size of box norms
(\ref{gvnbound})
\begin{equation*}
\abs{ \delta \operatorname T (f,S,f)}
\le \upsilon '+\delta (\delta _X \delta _Y \delta _D) ^2
\norm f.\Box.\norm f.\Box,X.
\le \kappa \delta ^{7/2}(\delta _X \delta _Y \delta _D) ^2\,.
\end{equation*}
The other two terms admit a similar bound.

We bound the final term  $ \operatorname T (f,f,f)$, with the inequality
Lemma~\ref{l.T}:
\begin{equation}\label{e.3f}
\Abs{\operatorname T (f,f,f)}\le   \upsilon '+
\delta ^{1/2}(\delta _{X}   \delta _Y   \delta _D) ^2
\norm f.\Box,X. \norm f.\Box,Y.   \le \kappa (\delta _X \delta _Y \delta _D) ^2 \delta ^{3}\,.
\end{equation}
Using the hypothesis on the Box norm, (\ref{gvnbound}),
and the equation (\ref{e.1f})
will prove the Lemma.

\end{proof}

\begin{lemma}\label{l.T}  For $ f_j\in \{f, S\}$ we have the estimate
\begin{equation}\label{e.T}
\abs{ \operatorname T (f_0,f_1,f_2)}
\le \upsilon '+ (\delta _X \delta _Y \delta _D) ^2\cdot
\begin{cases}
\norm f_0.2. \cdot \norm f_1.\Box,Y. \cdot \norm f_2.\Box,X.
\\
\norm f_0.\Box. \cdot \norm f_1.\Box,Y. \cdot \norm f_2.2.
\\
\norm f_0.\Box. \cdot \norm f_1.2. \cdot \norm f_2.\Box,X.
\end{cases}
\end{equation}
Here, $ \norm f.2.=\delta ^{1/2}$ while $ \norm S.2.=1$.
\end{lemma}

\begin{proof}

We prove an instance of the claimed inequalities:
\begin{equation}\label{e.instance}
\abs{ \operatorname T (f_0,f_1,f_2)}\le \upsilon '+
\norm f_0.2. (\delta _X \delta _Y \delta _D) ^2
\norm f_1.\Box,Y. \cdot \norm f_2.\Box,X.\,.
\end{equation}
By a change of  basis, this inequality implies the other
two.

Apply Cauchy Schwartz once, in the variables $(x,y)$, to get
\begin{align*}
\abs{ \operatorname T(f_0,f_1,f_2)} &\le \bigl(\mathbb E _{x,y}
f_0(x,y)^2\bigr) ^{1/2}\,\cdot\, \operatorname U ^{1/2}
\end{align*}
The first term on the right is no more than $  \norm f_0.2.( \delta _X  \delta _Y\delta _D)
^{1/2}$.
 As for the second term, it is
\begin{align*}
\operatorname U& \eqdef    \mathbb E _{x,y} D(x+y)
\Abs{ \mathbb E _{s}f_1(y+s,y) f_2(x,x+s)}^2
\\
&= \mathbb E _{y,s,s'}  \bigl\{\mathbb E _{x} D (x+y) f_1(y+s,y)f_1(y+s',y)
\} \cdot \{f_2(x,x+s) f_2(x,x+s') \}
\end{align*}

Note that in the definition of $\operatorname U$, we have inserted the
term $D(x+y)$, which arises from $A(x,y)$.\footnote{Without this term, we
would not get the right power of $ \delta _D$ in our estimates.}
We apply Proposition~\ref{uniformconvolution}
to  replace  $D(x+y)$  by $ \delta _D$ and then use Cauchy Schwartz again in
the variables $ x,x'$ to get
\begin{align*}
\operatorname U&=\delta_D
(\operatorname U_1\cdot \operatorname U_2) ^{1/2}+O ( \upsilon' );
\\
\operatorname U_1& \eqdef  \mathbb E _{y,y'\in H}
\Abs{\mathbb E _{s\in H}f_1(y+s,y)f_1(y'+s,y') } ^2
= \delta _X ^4 \delta _Y ^{2}\delta _D ^2   \norm
f_1. \Box,Y. ^{4}\,;
\\
\operatorname U_2& \eqdef \mathbb E _{x,x'\in H} \abs{ \mathbb E _{s\in H}
f_2(x,x+s)f_2(x',x'+s)}^2
=\delta _X ^2 \delta _Y ^{4}\delta _D ^2   \norm
f_2. \Box,X. ^{4}\,.
\end{align*}
A change of variables permits the identification of $ U_1 $ and $U_2 $.
This completes the proof of (\ref{e.3f}).

\end{proof}

\section{Proof of Density Increment Lemma}

\begin{proof}[Proof of Lemma~\ref{increment}.]

We
prove this assertion:  Fix $ c <1$.  There is a 
constant $ \kappa  $  so that the following holds.
Assume that $ \norm f.\Box.\ge c\delta ^{5/4}$,
and show that there are subsets $ X'\subset X$, $ Y'\subset Y$ so that
\begin{gather*}
\mathbb P (X'\mid X)\,,\, \mathbb P (Y'\mid Y)\ge \kappa \delta ^{2 }\,,
\\
\mathbb P (A\mid X' \times Y' \cap X'\stackrel {\textup{diag}} \times  D)
\ge \delta +\kappa  \delta ^{2 }\,.
\end{gather*}
We emphasize that the box norm we use is the one given in (\ref{e.box}).
And we will \emph{not} refine the diagonal coordinate.
This is one instance of the Lemma, which
by a change of coordinates,
this will prove the Lemma as stated.

\smallskip

We can assume that the fibers above points $ x\in X$, and $ y\in Y$,
typically behave as expected.
Namely, we assume that
\begin{equation}\label{e.8}
\begin{split}
\mathbb E _{x\in X} \Abs{\mathbb E _{y\in Y} f (x,y)} ^2 &
\le \kappa ^3  \delta ^{2} \delta _D ^2 \,,
\\
\mathbb E _{y\in Y} \Abs{\mathbb E _{x\in X} f (x,y)} ^2 &\le \kappa ^3  \delta ^{2} \delta _D ^2 \,.
\end{split}
\end{equation}
for otherwise, we can apply Lemma~\ref{l.fibers} to conclude the Lemma.

For a point $ (x,y)\in A$, consider $ N_x \eqdef \{y'\mid (x,y')\in A\}$,
$ N_y \eqdef \{x'\mid (x',y)\in A\}$.  These are the neighbors of $ x$ and of $ y$, respectively.
We need to consider points for which
these sets are about as big as they should be.
Set
\begin{align*}
X'' &= \{x\in X \mid
  \mathbb P (N_x|Y) \ge \kappa  \delta ^{5}  \delta _D  \,,\,
  \Abs{\mathbb E _{y\in Y} f (x,y) } < \kappa \delta \delta _D  \}
 \,,
\\
Y''&= \{y\in Y \mid
  \mathbb P (N_y|X) \ge \kappa  \delta ^{5}  \delta _D    \,,\,
    \Abs{\mathbb E _{x\in X} f (x,y) } < \kappa \delta \delta _D  \}
  \,.
\end{align*}
It is clear that  these sets are most of $ X$ and $ Y$ respectively.
In particular, in view of (\ref{e.8}) we have
\begin{equation}\label{e.88}
\mathbb P (X''\mid X)\,,\, \mathbb P (Y''\mid Y)\ge
1- \delta ^{2} \,.
\end{equation}
Clearly, we can assume that
\begin{equation} \label{e.7}
 \mathbb E _{\substack{x\in X''\\ y\in Y'' }} f (x,y)  <\kappa \delta ^{2} \delta _D
\end{equation}
for otherwise we already proved the Lemma.  But, it is also the case that we can
assume
\begin{equation}\label{e.888}
-\kappa \delta ^{4}\delta _D  <\mathbb E _{\substack{x\in X''\\ y\in Y'' }} f (x,y)\,
\end{equation}
for if this inequality fails, we apply Lemma~\ref{l.lowDensity} to conclude the proof
of the Lemma.

We further note  that we have
\begin{equation} \label{e.P}
\mathbb E _{\substack{x,x'\in X''\\ y,y'\in Y'' }}
f (x,y) f (x',y) f (x,y') f (x',y')\ge \tfrac {c ^{4}} 2 \delta ^{5}\,.
\end{equation}
If we were taking the expectation over $ X$ and $ Y$, this would
follow from the assumption that $ \norm f.\Box.\ge c \delta ^{5/4}$.
As we are taking the expectation over $ X''$ and $ Y''$,
we need to show that taking the expectation over the  complement of $ X''$
we get an appropriate upper bound, namely
\begin{align} \label{e.Complement}
\Abs{ \mathbb E _{\substack{x,x'\in X\\ y,y'\in Y'' }} \mathbf 1 _{X-X''} (x)
f (x,y) & f (x',y) f (x,y') f (x',y')}
\le \upsilon '+ \kappa \delta ^{5} \delta _D ^{4}\,.
\end{align}
Three  similar inequalities hold, so using the assumption that
$ \norm f.\Box.\ge c \delta ^{5/4}$,
taking $ 0<\kappa <\tfrac {c ^{4}} 8 $ we will see that  (\ref{e.P}) holds.

To see (\ref{e.Complement}), first  observe that by definition of $ X''$,
\begin{equation}\label{e.Comp1}
\mathbb E _{\substack{x , x'\in X\\ y,y'\in Y'' }} \mathbf 1 _{X-X''} (x)
N_x (y) D (x+y') D (x'+y) D (x'+y')
\le
 \upsilon '+ \kappa \delta ^{5} \delta _D ^{4}\,.
\end{equation}
This does not prove (\ref{e.Complement}) since $ f (x,y)$ is not supported on
$ N_x (y)$.  But, we also have the similar inequality
\begin{equation}\label{e.Comp2}
\mathbb E _{\substack{x , x'\in X\\ y,y'\in Y'' }} \mathbf 1 _{X-X''} (x)
 D (x+y) N_x (y')  D (x'+y) D (x'+y')
 \upsilon '+ \kappa \delta ^{5} \delta _D ^{4}\,.
\end{equation}
Of course $ f=\delta S - A$, so we can estimate
\begin{align*}
\delta \Abs{\mathbb E _{\substack{x , x'\in X\\ y,y'\in Y'' }}  & \mathbf 1 _{X-X''} (x)
D (x+y) f (x,y') f (x',y) f (x',y') }
\\
& \le
\upsilon '+ \kappa \delta ^{5} \delta _D ^{4}
+ \delta ^2  \Abs{\mathbb E _{\substack{x , x'\in X\\ y,y'\in Y'' }} \mathbf 1 _{X-X''} (x)
D (x+y)  D(x+y') f (x',y) f (x',y') }
\\
& \le
\upsilon '+ \kappa \delta ^{5} \delta _D ^{4}
+ \delta ^2 \delta ^2 _D   \Abs{\mathbb E _{\substack{  x , x'\in X\\ y,y'\in Y'' }}
\mathbf 1 _{X-X''} (x)
f (x',y) f (x',y') }
\\
& \le
\upsilon '+ \kappa \delta ^{5} \delta _D ^{4}
+ \delta ^4 \delta ^2 _D   \mathbb E _{\substack{  x'\in X }}
\Abs{
 \mathbb E _{y\in Y''} f (x',y)  } ^2
 \\
 & \le \upsilon '+ \kappa \delta ^{5} \delta _D ^{4}
+ \kappa \delta ^6 \delta ^4 _D\,.
\end{align*}
This completes the proof of (\ref{e.Complement}).

\medskip

To find the subset on which $ A$ has increased density, we
consider any subset of the form $  N_y \times  N_x\cap  N_y
\stackrel {\textup{diag}} {\times } D$,
where $ y\in Y''$ and $ x\in X''$. 
Clearly we are interested in the
largest increase in density, for which we
estimate
\begin{gather} \label{e.inf}
\sup _{\substack{x\in X'', y\in Y'' \\ (x,y)\in A}}
\frac {\abs{ A \cap \{N_y \times N_x \cap S''\} }}
	{\abs{ N_y \times N_x \cap S''}}
\ge \frac {\operatorname Q'' (A,A,A,A)} {\operatorname Q'' (A,A,A,S)} \,,
\\ \nonumber
S''  \eqdef
X'' \times Y'' \cap X'' \stackrel {\textup{diag}}\times D\,,
\\ \nonumber
\operatorname Q'' (f_0,f_1,f_2,f_3)
 \eqdef
\mathbb E _{ \substack{x,x'\in X'' \\ y,y''\in Y''}}
f_0 (x,y)f_1 (x',y)f_2 (x,y')f_3 (x',y')\,.
\end{gather}
Here, we note that for $ (x,y)\in X'' \times Y'' $ we have 
$\abs{ A \cap \{N_y \times N_x \cap S''\} }>0 $, so that we 
are not dividing by zero in (\ref{e.inf}).  By definition, we have 
\begin{equation*}
\mathbb P (N_y \times N_x\mid X \times Y)\ge k ^2 \delta ^{10} \delta _D ^2 \,.
\end{equation*}
And by uniformity, we have 
\begin{align*}
\mathbb P ( N_y \times N_x \cap S'')
&= \mathbb E _{x',y'\in H} N _{y} (x') N _{x} (y') D (x'+y')
\\&\ge \upsilon '+ k ^2 \delta ^{10} \delta _D ^{3} \delta _X \delta _Y>0\,. 
\end{align*}
In the last inequality, $ \upsilon '$ is a function of the uniformity constant.

There is a gain of regularity in passing to the expectations on
the right hand side of (\ref{e.inf}).
Expand $ A=f+\delta S$ in the last place in $ \operatorname Q (A,A,A,A)$ to see that
\begin{align*}
\frac {\operatorname Q'' (A,A,A,A)} {\operatorname Q'' (A,A,A,S)}
\ge \delta +\frac {\operatorname Q'' (A,A,A,f)} {\operatorname Q'' (A,A,A,S)}
\end{align*}
and so we should show that the last fraction is at least $ c '
\delta ^2 $,  which we do by showing that
\begin{gather}\label{e.Qtop}
\operatorname Q'' (A,A,A,f)\ge c' \delta ^{5} \delta _D ^{4}\,,
\\
\label{e.Qbottom}
0<\operatorname Q'' (A,A,A,S) < 10 \delta ^{3} \delta _{D} ^{4}\,.
\end{gather}
This we will do, assuming one more condition.  If this last condition fails, we will
get a density increment of $ \delta ^{2 }$.

Define
\begin{equation*}
\alpha _{X} ^2 \eqdef \delta  \delta _D ^2  \mathbb E _{x\in X''} \Abs{\mathbb E _{y \in Y''} f (x,y)} ^2 \,,
\qquad
\alpha _{Y} ^2 \eqdef \delta   \delta _D ^2 \mathbb E _{y\in Y''} \Abs{\mathbb E _{x\in X''} f (x,y)}^2 \,,
\end{equation*}
By (\ref{e.8}) and the definitions of $ X''$ and $ Y''$, one can see that
these two quantities are at most $ \kappa \delta ^{3} \delta ^{4}_D$, and hence are only
a small fraction of the major term in the considerations below.

Using $ A=\delta S+f$,  $ \operatorname Q'' (A,A,A,S)$ has the expansion
\begin{align} \label{e.w1}
\operatorname Q'' (A,A,A,S)&= \delta ^{3} \mathbb E _{\substack{x,x'\in X'' \\ y,y'\in Y''}}
D (x+y)D (x+y')D (x'+y)D (x'+y')				
\\ \label{e.w2}
&\quad + 3 \delta ^{2} \mathbb E _{\substack{x,x'\in X'' \\ y,y'\in Y''}}
f (x,y)D (x+y')D (x'+y)D (x'+y')
\\ \label{e.w3}
&\quad + \delta
\mathbb E _{\substack{x,x'\in X'' \\ y,y'\in Y''}}
f (x,y)f (x,y')D (x'+y)D (x'+y')
\\ \label{e.w4}
&\quad + \delta
\mathbb E _{\substack{x,x'\in X'' \\ y,y'\in Y''}}
D (x+y)f (x,y')f (x',y)D (x'+y')
\\ \label{e.w5}
&\quad + \delta
\mathbb E _{\substack{x,x'\in X'' \\ y,y'\in Y''}}
f(x,y)D (x+y')f (x',y)D (x'+y')
\\ \label{e.w6}
&\quad +
\mathbb E _{\substack{x,x'\in X'' \\ y,y'\in Y''}}
f (x,y)f (x,y')f (x',y)D (x'+y')\,.
\end{align}
Clearly, the uniformity in $ D$ is relevant.
The right hand side of (\ref{e.w1}) is $ \delta ^{3} \delta _D ^{4}$ plus
a term controlled by uniformity;
the term in (\ref{e.w2}) is, by (\ref{e.888}), at least $-3 \kappa {\delta ^{6}}$,
plus a term controlled by uniformity;
(\ref{e.w3}) is $  \alpha _{X} ^2 $, plus a term controlled by uniformity;
(\ref{e.w4}), by (\ref{e.7}) and (\ref{e.888}), obeys the inequality
\begin{equation*}
(\ref{e.w4})< \upsilon '+\delta \delta _D ^2 \Abs{\mathbb E _{x\in X,y\in Y} f (x,y)} ^2< \upsilon ' + \kappa ^{2}\delta ^{3} \delta _D ^{4}\,;
\end{equation*}
(\ref{e.w5}) is approximately in $  \alpha _Y ^2 $;
while the last term (\ref{e.w6}) is not one that admits an obvious control, and we write
\begin{equation} \label{e.DDD}
(\ref{e.w6})=\upsilon '+  \Delta \,, \qquad
\Delta  \eqdef  \delta _D  \mathbb E _{\substack{x,x'\in X'' \\ y,y'\in Y''}}
f (x,y)f (x,y')f (x',y)\,.
\end{equation}
We can assume that $ \abs{ \Delta }\le \kappa \delta ^{3}  \delta _D ^4$, otherwise
we apply Lemma~\ref{l.Delta} to finish the proof of the Lemma, getting
a density increment of the order of $ \delta ^{2 }$.
This proves (\ref{e.Qbottom}).

\smallskip

The expression $ \operatorname Q (A,A,A,f)$ admits a very similar expansion.
\begin{align} \label{e.W1}
\operatorname Q'' (A,A,A,f)&= \delta ^{3} \mathbb E _{\substack{x,x'\in X'' \\ y,y'\in Y''}}
D (x+y)D (x+y')D (x'+y)f(x',y')				
\\
& \quad +\operatorname  Q_2 + \operatorname Q_3
\\
 \label{e.W6}
&\quad +
\mathbb E _{\substack{x,x'\in X'' \\ y,y'\in Y''}}
f (x,y)f (x,y')f (x',y)f(x',y') \, ,
\\ \label{e.W2}
 \operatorname  Q_2 & \eqdef  \delta ^{2}\bigl\{
 \mathbb E _{\substack{x,x'\in X'' \\ y,y'\in Y''}}  f (x,y)D (x+y')D (x'+y)f(x',y')
\\ \nonumber
&\quad +
\mathbb E _{\substack{x,x'\in X'' \\ y,y'\in Y''}}
D(x+y)f (x,y')D (x'+y)f(x',y')
\\ \nonumber
&\quad +
\mathbb E _{\substack{x,x'\in X'' \\ y,y'\in Y''}}
D (x+y)D (x+y')f (x',y)f(x',y') \bigr\} \,,
\\ \label{e.W5a}
 \operatorname Q_3 & \eqdef  \delta \{
\mathbb E _{\substack{x,x'\in X'' \\ y,y'\in Y''}}
D(x+y)f(x,y')f (x',y)f(x',y')
\\ \nonumber
&\quad +
\mathbb E _{\substack{x,x'\in X'' \\ y,y'\in Y''}}
f(x,y)D (x+y')f (x',y)f(x',y')
\\ \nonumber
&\quad +
\mathbb E _{\substack{x,x'\in X'' \\ y,y'\in Y''}}
D(x+y)f (x,y')f (x',y)f(x',y') \} \,.
\end{align}
The right hand side of (\ref{e.W1}) is greater than $- \kappa \delta ^{3} \delta_D ^{4}$, plus  a term controlled by uniformity by (\ref{e.7}) and (\ref{e.888});
the term in (\ref{e.W6}) is the box norm, which by (\ref{e.P}) is
at least $ c' \delta ^{5} \delta _D ^{4}$; here of course, $ c'>0$ is fixed in
advance, and we take  $ \kappa $ much smaller than $ c'$; the terms which make up the definition
of $ \operatorname Q_2$ all involve two $ f$'s, and after taking uniformity into
account each individual term is  positive, and
so can be ignored as we are obtaining a lower bound for
$ \operatorname Q (A,A,A,f)$; the three terms in the
definition of $ \operatorname Q_3$ are all of the form $ \upsilon '+\Delta $
where $ \Delta $ is defined in (\ref{e.DDD}).  In particular, we
have already assumed $ \abs \Delta \le \kappa \delta ^{3} \delta _D ^{4}$.
For $ 0<\kappa $ sufficiently small, this proves (\ref{e.Qtop}), and so the proof
 of this Lemma.

\end{proof}

We use the following simple
variant of the Paley Zygmund inequality.  It states, in particular, that a random variable,
bounded in $ L ^{\infty }$ norm by one, with mean zero, and standard deviation $ \sigma $,
must be at least a constant multiple of $ \sigma $ on a set of probability proportional
to the variance $ \sigma ^2 $.

\begin{proposition}\label{variance} Let $ 1<p<\infty $.  Then there is a $ c>0$ 
so that for all $ 1\le p <\infty $, and random variables 
  $ Z$  with $ -1\le Z\le 1$,  $ \mathbb E
Z=0$, and $ \mathbb E \abs Z^p= \sigma ^p$. Then, $ \mathbb P (Z>c \sigma^p
)\ge c_p \sigma ^p$.
\end{proposition}


\begin{proof}
In fact, we can take $ c=\tfrac 15$.  We assume that the conclusion is 
false and seek a contradiction. 
Since $ \mathbb E Z=0$ we have 
\begin{align*}
- \mathbb E Z \mathbf 1 _{ \{Z<0\}} & = 
\mathbb E Z \mathbf 1 _{ \{Z>0\}}
\\
& \le \mathbb P (Z> c \sigma ^{p})+ \mathbb E Z  \mathbf 1 _{ 
\{0<Z< c \sigma ^{p}\}}
\\
& \le 2 c \sigma ^{p}\,. 
\end{align*}
With this, and the  fact that $ p\ge 1$ while $ Z$ is bounded by $ 1$, 
we can now estimate 
\begin{align*}
\sigma ^{p}= \mathbb E \lvert  Z\rvert ^{p} & =
 \mathbb E \lvert  Z\rvert ^{p} \mathbf 1 _{ \{Z<0\}} 
 +  \mathbb E \lvert  Z\rvert ^{p} \mathbf 1 _{ \{Z>0\}} 
\\
& \le 2  \mathbb E Z \mathbf 1 _{ \{Z>0\}} \le 4 c \sigma ^{p}\,. 
\end{align*}
This is a contradiction. 

\end{proof}

\begin{lemma}\label{l.fibers} If it is the case that
\begin{equation*}
\mathbb E _{x\in X} \Abs{\mathbb E _{y\in Y} f (x,y)} ^2 \ge c \delta ^{2} \delta _D ^2
\end{equation*}
then there is a set $ X'\subset X $ for which $ \mathbb P (X'\mid X)\ge 
\tfrac c  {12}\delta ^{2}$
on which we have
\begin{equation*}
\mathbb P (A\mid X' \times Y\cap X' \stackrel {\textup{diag}} \times D)
\ge \delta +\tfrac c {12} \delta ^2 \,.
\end{equation*}
\end{lemma}

\begin{proof}
The set of $ x\in X$ for which
\begin{equation*}
2\delta _D 
\le
\Abs{\mathbb E _{y\in Y} f (x,y)} 
\le
 \Abs{\mathbb E _{y\in Y} D (x+y)} 
\end{equation*}
has probability that is controlled by uniformity, and hence is negligible.
Thus, the Lemma follows immediately from the Paley Zygmund inequality.
\end{proof}

\begin{lemma}\label{l.lowDensity} Suppose that there is a
sublattice $ S''=X'' \times Y''  \cap X'' \stackrel {\textup{diag}} \times D  $ with
\begin{align}\label{e.low1}
\mathbb P (X''\mid X )\,,\, \mathbb P (Y''\mid Y)&\ge 1-\lambda\,,
\\ \label{e.low2}
\mathbb P (A\mid S'')&\le \delta -\tau \,.
\end{align}
Then, there is a sublattice $ S'=X' \times Y' \cap X' \stackrel {\textup{diag}} \times D$
on which
\begin{align}\label{e.low3}
\mathbb P (X'\mid X)\,,\, \mathbb P (Y'\mid Y)&\ge \tfrac12 \tau \,,
\\ \label{e.low4}
\mathbb P (A\mid S')&\ge \delta + \kappa \tau \lambda ^{-1}  \,.
\end{align}
\end{lemma}


\begin{proof}
Notice that the
density of $ A$ on $ S-S''$ must be strictly larger than $ \delta $. Namely,
\begin{align*}
\mathbb P (A\mid S-S'')&\ge \frac{\delta[1- \mathbb P (S''\mid S)]
+\tau \mathbb P (S''\mid S)}
{1-\mathbb P (S''\mid S)}
\\&\ge \delta + \tau \frac{\mathbb P (S''\mid S)} {1-\mathbb P (S''\mid S)}
\\&\ge \delta +\kappa \tau \lambda ^{-1}\,.
\end{align*}
But, clearly $ S-S''$ is a union of three sublattices, and on one of these three,
$ A$ must have density at least $ \delta + \kappa \tau \lambda ^{-1} $.
We must have $ \mathbb P (X''\mid X)\le
1-\tfrac12 \tau $, otherwise we contradict (\ref{e.low2}). This finishes the proof.
\end{proof}


\begin{lemma}\label{l.Delta}  Fix $ c>0$.  If it is the case that
\begin{equation}\label{e.Delta}
\Abs{\mathbb E _{x,y} f (x,y) \mathbb E _{x'} f (x',y) \, \mathbb E _{y'} f (x,y')}
\ge c \delta ^{3} \delta _D ^{3}
\end{equation}
Then, there is a constant  $ c'=c' (c)$, $ X'\subset X$, $ Y'\subset Y$ with
\begin{gather}\label{e.Delta1}
\mathbb P (X'\mid X)\,,\, \mathbb P (Y'\mid Y)\ge c'\delta ^{2 }\,,
\\
\mathbb P (A\mid X' \times Y'\cap X' \stackrel {\textup{diag}} \times  D )
\ge \delta +  c'\delta ^{2 }\,.
\end{gather}
\end{lemma}

\begin{proof}
By uniformity in $ D$, we have
\begin{equation*}
\mathbb P \bigl( \Abs{\mathbb E _{y} f (x,y)} > 2 \delta _D \bigr)
\le \upsilon '\,.
\end{equation*}
That is, the effective $ L ^{\infty }$ bound on $ \mathbb E _{y} f (x,y)$
is $ 2 \delta _D$.

H\"olders inequality and the assumption (\ref{e.Delta})
implies
\begin{align*}
c \delta ^{3} \delta _D ^{3}&\le \bigl[ \mathbb E _{x,y} \lvert  f (x,y)\rvert  ^{3} \bigr] ^{1/3}\cdot
\bigl[  \mathbb E _{x,y} D (x+y) \Abs{\mathbb E _{y} f (x,y)} ^{3/2}
\cdot \Abs{\mathbb E _{x} f (x,y)} ^{3/2}  \bigr] ^{2/3}
\\
& \le 2 \delta ^{1/3 } \delta _D  ^{1/3}
\bigl[ \mathbb E _{x,y} D (x+y) \Abs{\mathbb E _{y'} f (x,y')} ^{3/2}
\cdot \Abs{\mathbb E _{x'} f (x',y)} ^{3/2}  \bigr] ^{2/3}
\\
&\le \upsilon '+ 2 \delta ^{1/3 } \delta _D
\bigl[ \mathbb E _{x}\Abs{\mathbb E _{y'} f (x,y')} ^{3/2}
\cdot \mathbb E _{y}\Abs{\mathbb E _{x'} f (x',y)} ^{3/2}  \bigr] ^{2/3}\,.
\end{align*}
Here it is essential that we insert the term $ D (x+y)$ when we apply
H\"olders inequality.  Note that uniformity in $ D$ is then used to obtain the
full power of $ \delta _D$ in the last line.

Thus, we must have e.\thinspace g.~
\begin{equation*}
 \mathbb E _{y} \Abs{\mathbb E _{x} f (x,y)} ^{3/2}  
\ge c''
\delta ^{2 } \delta _D ^{3/2}
\end{equation*}
Then, the conclusion of our Lemma follows from the Paley Zygmund inequality.

\end{proof}

\section{Proof of Lemma~\ref{uniformity}} 

We include a proof of this Lemma as we are requiring a  more
than is claimed in e.\thinspace g.~Ben Green's survey \cite {MR2187732}.
In particular, we claim that all three sets $ X,Y,D$ can be
uniformized.

Let us indicate the central way that uniformity is used in this proof.
See \cite[Lemma 3.4(1)]{MR2187732}.

\begin{proposition}\label{p.uniUp}  For any subset $ X\subset H$,
there is a partition of $ H$ into two affine subspaces $ H'$ and $ H''$
for which
\begin{equation}\label{e.uniUp}
\tfrac 12 \bigl\{
\mathbb P (X\mid H')^2+  \mathbb P (X\mid H'') ^2
\bigr\}
\ge \mathbb P (X\mid H) ^2 + \tfrac 18\norm X.\textup{Uni}. ^2 \,.
\end{equation}
\end{proposition}

The following more technical Lemma describes
a key inductive procedure in the proof.

\begin{lemma}\label{l.X} Let $ 0<t,u<1$ be positive parameters.
Let $ U\subset H$. Suppose that $ \operatorname {dim}(H)\ge 10 (  t u ^2  ) ^{-1} $.
Then, there is a partition $ \mathcal P$ of $ H$,
so that writing $ \mathcal P=\mathcal U\cup \mathcal N$ ,
\begin{enumerate}
\item   All $ H'\in \mathcal P$ have  dimension
$\operatorname {dim} (H')\ge \operatorname {dim}(H)-4(  t u ^2  ) ^{-1}$;
\item  $ \norm U. \textup{Uni}.\le u$ for all $ H'\in \mathcal U$;
\item For $ \delta _U= \mathbb P (U\mid H)$, 
\begin{equation*}
\mathbb P \bigl( \bigcup \{H'\mid H'\in \mathcal N\} \mid H\bigr)\le t \delta _U\,.
\end{equation*}
\end{enumerate}
\end{lemma}


\begin{remark}\label{r.tu}  We will not use this Lemma as stated,
but rather the more complicated variant that follows.  We include
this statement and proof for clarity's sake.
Some of the notation above is taken from Ben Green's note \cite {MR2187732}.
In particular, $ \mathcal U$ is for `uniform', and
$ \mathcal N$ is for `non--uniform.'
\end{remark}

\begin{proof}
The proof is an inductive procedure, though we will not define
the collections of hyperspaces $ \mathcal U$ and $ \mathcal N$ until the
conclusion of the iteration. Initialize variables
\begin{equation*}
\mathcal Q\leftarrow \{H\}.
\end{equation*}
Also initialize a counter $ m\leftarrow 0$.
Given $ \mathcal Q$, define $ \mathcal R$ to be those $ H'\in \mathcal Q$
for which $\norm U\cap H'. \textup{Uni}.\ge u$.
\texttt {WHILE}
\begin{equation*}
 \mathbb P \bigl( \bigcup \{H' \mid H'\in  \mathcal R \} \bigr)\ge t \delta _U
\end{equation*}
update $ m\leftarrow m+1$.  And for each $ H'\in \mathcal R$
  apply Lemma~\ref{l.X} to $ H'$.  Thus, $ H'=H'_1\cup H'_2$ for which
\begin{equation}\label{e.Up}
\tfrac 12 \bigl\{
\mathbb P (U\mid H_1')^2+  \mathbb P (U\mid H'_2) ^2
\bigr\}
\ge \mathbb P (U\mid H') ^2 + \tfrac 14 u ^2 \,.
\end{equation}
Update $ \mathcal Q\leftarrow \mathcal Q- \mathcal R\cup \{H'_1,H'_2 \mid H'\in \mathcal R\}$.
When each $ H'\in \mathcal R$ has been so split, set $ \mathcal P_m\leftarrow \mathcal Q$.
The \textsf{WHILE} loop then repeats.

Once the \textsf{WHILE} loop stops, return the value of $ m$,
the sequence of partitions $ \mathcal P_1,\dotsc, \mathcal P_m$, and
 $ \mathcal Q=\mathcal P_m$.
Define
\begin{align*}
\mathcal U &\eqdef \{H'\in \mathcal Q\mid \norm X\cap H'.\textup{Uni}.\le u\}\,,
\end{align*}
and $ \mathcal N \eqdef \mathcal Q-\mathcal U$.
The subspaces in $ \mathcal Q$ have dimension at least
$ \operatorname {dim} (H)-m$.

\medskip

Once the \texttt {WHILE} loop stops, it is clear that the conclusions
$ (2)$ and $ (3)$ of the Lemma hold.  We concern ourselves with the first
conclusion about the dimension of the subspaces involved.  Observe that
at each iterate of the \texttt {WHILE} loop, the dimensions of the partition
can only lose one dimension.  It suffices to estimate the number of iterates
of the Loop, which is the counter $ m$.

Consider the quantities
\begin{equation*}
S_j=\sum _{H'\in \mathcal P_j} \mathbb P (U\mid H') ^2 \mathbb P (H'\mid H)
=\mathbb E \; \mathbb E (U\mid \mathcal P_j) ^2 \,,
\end{equation*}
which are the mean square densities of $ U$ relative to the partitions $ \mathcal P_j$.
(Here we  are relying on the usual notations for conditional second moments.)
Obviously these quantities are  less than $ \delta _{U}=\mathbb P (U\mid H)$.
Each $ H'$ which is split in (\ref{e.Up}), the mean square density in $ H'$ is
increased by $ \tfrac 14 u ^2  $.
And, this takes place, at each iterate of the loop, in
a portion of the whole space $ H$ that is at least probability $ t \delta _U$.
From this, we see that $ S _{j+1}\ge S_j+\tfrac14 t \delta _U u ^2 $.  Consequently,
$ m t \delta _U u ^2 \le S_m\le \delta _U $.  This proves $ (1)$, and
so the proof is complete.

\end{proof}


The more technical statement that we need is as follows.  Whereas in the
first Lemma, we have a single subset $ U\subset H$ and construct
a `good' partition, in this statement we have three subsets $ U_1,U_2, U_3 \subset H$
and construct a single `good' partition in the product space $ H\times H$.

\begin{lemma}\label{l.UU} Let $ 0<t,u<1$ be positive parameters.
Let $ U_1, U_2, U_3\subset H$. Suppose that $ \operatorname {dim}(H)\ge 10 (
t u ^2  ) ^{-1} $.
Then, there is a partition $ \mathcal P$ of $ H\times H$,
so that writing $ \mathcal P= \mathcal U\cup \mathcal N_1\cup \mathcal N_2\cup \mathcal N_3$,
we have the following.
\begin{enumerate}
\item   For all $ V_1\times V_2\in \mathcal P$,
\begin{equation*}
\operatorname {dim} (V_1)=\operatorname {dim} (V_2)\ge \operatorname {dim}(H)-
4(    t u ^2  ) ^{-1}\,.
\end{equation*}
\item  Moreover, $ V_1,V_2$ as affine subspaces of $ H$, are translates of each other.
\item For all $ V_1\times V_2\in \mathcal U$, $\max _{j=1,2} \norm U_j\cap V_j. \textup{Uni}.\le u$
and $ \norm U_3\cap (V_1+V_2). \textup{Uni}.\le u $.
\item
\begin{gather*}
\mathbb P \bigl( \bigcup \{V_j\mid V_1\times V_2
\in \mathcal N_j\} \mid H\bigr)\le t \delta _{U_j}\,\ j=1,2
\,,
\\
\mathbb P \bigl( \bigcup \{V_1+V_2\mid V_1\times V_2
\in \mathcal N_3\} \mid H\bigr)\le t \delta _{U_3}\,.
\end{gather*}
\end{enumerate}
\end{lemma}

\begin{remark}\label{r.UU}
If we did not insist on the second conclusion, that the relevant
subspaces be translates of one another, we could simply take a product of
partitions arising from Lemma~\ref{l.X}.  Due to the prominent role of the
diagonals in our problem, this is an essential point for us.
\end{remark}

\begin{proof}

We describe the recursive procedure.
Initialize variables
\begin{equation*}
\mathcal Q\leftarrow \{H\times H\}.
\end{equation*}
and three counters $ m_j\leftarrow 0$, $ j=1,2,3$.
The recursive procedure we describe will return the partition we
need, as well as some auxiliary data that we need to prove the Lemma.

Given $ \mathcal Q$, define $ \mathcal R_j$, $ j=1,2,3$  by
\begin{gather*}
\mathcal R _{j} = \{V_1 \times V_2 \in \mathcal Q \mid \norm U_j\cap V_j.
\textup{Uni}.\ge u\}\,, \qquad j=1,2\,
\\
\mathcal R _{3} = \{V_1 \times V_2 \in \mathcal Q \mid \norm U_3\cap (V_1+V_2).
\textup{Uni}.\ge u\}\,.
\end{gather*}

For $ j=1,2$
\begin{equation} \label{e.Nj}
\textsf {IF}\qquad \mathbb P \bigl( \bigcup\{ V_j \mid V_1
\times V_2 \in \mathcal R _{j}\} \bigr)\ge t \delta _{U_j}
\end{equation}
\textsf{THEN} update $ m_j\leftarrow m_j+1$.
\textsf {WHILE} $ \mathcal R_j\neq \emptyset$
\begin{itemize}
\item Apply Lemma~\ref{l.X} to $ V_j$.  Thus, $ V_j=V_j'\cup V_j''$ for which
\begin{equation}\label{e.jUp}
\tfrac 12 \bigl\{
\mathbb P ({U_j}\mid V_j')^2+  \mathbb P ({U_j}\mid V_j'') ^2
\bigr\}
\ge \mathbb P ({U_j}\mid V_j) ^2 + \norm {U_j}\cap V_j.\textup{Uni}. ^2 \,.
\end{equation}
\item (This point is a departure from the previous proof.)  Let $ k$ be the other
index.\footnote{That is, if $ j=1$, then $ k=2$. Notice that we are enforcing a
partition on $ V_k$ that does not necessarily have anything to do with
increasing conditional variances.}
Since $ V_k$ is a translate of $ V_j$, we can choose translates $ V_k'$
and $ V_k''$ of $ V_j'$ which also partition $ V_k$.
\item
At this point, let us observe that we will have $ V_1'+V_2',\,
V_1'+V_2'',\,  V_1''+V_2',\, V_1''+V_2''\subset V_1+V_2$.
And indeed, all of these subspaces are affine translates of $ V_1'$,
hence these four
subspaces are made of two equal pairs, which partition $ V_1+V_2$.
\item
Update
\begin{gather*}
 \mathcal Q\leftarrow (\mathcal Q- \{V_1\times V_2\})\cup \{V_1'\times V_2',
 V_1''\times V_2', V_1'\times V_2'', V_1''\times V_2''\}
\\
\mathcal R_j\leftarrow \mathcal R_j- \{V_1 \times V_2\}.
\end{gather*}
\end{itemize}


When the \textsf {WHILE} loop stops,
 define $ \mathcal P _{m_1+m_2+m_3} \eqdef  \mathcal Q$
and $ \sigma _{j} (m_j)=m_1+m_2+m_3$.

We now describe the procedure as applied to the set $ U_3$. That the products
in $ \mathcal Q$ are products of translates of the same subspace plays a critical
role in this formulation. In particular, for $V_1  \times V_2\in \mathcal P$,
$ V_1+V_2$ is a translate of $ V_1$ (and $ V_2$).
\begin{equation} \label{e.3Nj}
\textsf {IF}\qquad \mathbb P
\bigl( \bigcup \{V_1+V_2\mid  V_1 \times V_2 \in \mathcal R_3\}
 \bigr)\ge t \delta _{U_3}
\end{equation}
\textsf{THEN} update $ m_3\leftarrow m_3+1$.   For each $ V_1\times V_2\in \mathcal R_3$
\begin{itemize}
\item Apply Lemma~\ref{l.X} to $W=V_1+V_2$.
Thus, $W=W'\cup W''$ for which
\begin{equation}\label{e.3Up}
\tfrac 12 \bigl\{
\mathbb P ({U_3}\mid W')^2+  \mathbb P ({U_3}\mid W'') ^2
\bigr\}
\ge \mathbb P ({U_3}\mid W) ^2 + \norm {U_3}\cap W.\textup{Uni}. ^2 \,.
\end{equation}
\item
  Since the spaces $ V_j$ are translates of $ W$, we can choose translates $ V_j'$
and $ V_j''$ of $ W'$ which also partition $ V_j$, $ j=1,2$..
\item
Update
\begin{gather*}
 \mathcal Q\leftarrow (\mathcal Q- \{V_1\times V_2\})\cup \{V_1'\times V_2',
 V_1''\times V_2', V_1'\times V_2'', V_1''\times V_2''\}
\\
\mathcal R_3\leftarrow \mathcal R_3- \{V_1 \times V_2\}.
\end{gather*}
\end{itemize}
Repeat these steps until $ \mathcal R_3$ is exhausted.  Then,
define $ \mathcal P _{m_1+m_2+m_3} \eqdef \mathcal Q$
and $ \sigma _{3} (m_j)=m_1+m_2+m_3$.

Iteratively apply the three conditionals, two in (\ref{e.Nj}) and one in (\ref{e.3Nj}).
\textsf{STOP} when all three conditionals fail.
Return the values of $ m_j$, the sequence
of partitions of  $ \mathcal P _{j}$ for $ 1\le j\le m_1+m_2$,
the `partition times' $ \{\sigma _j (n)\mid 1\le n\le m_j\}$
and collection   $ \mathcal Q$.

Define
\begin{equation*}
\mathcal {N}_j \eqdef \{V_1\times V_2\in \mathcal Q\mid \norm U_j\cap V_j.\textup{Uni}.\ge u\}\,,
\qquad j=1,2
\end{equation*}
and similarly define $ \mathcal N_3$.
Set $ \mathcal U=\mathcal Q-\mathcal N_1-\mathcal N_2- \mathcal N_3$.

All subspaces chosen in this way have
dimension at least equal to $ \operatorname {dim} (H)-m_1-m_2-m_3$.
We need only provide upper bounds on the   $ m_j$, as all the other
claims of the Lemma are evident from the construction.

We claim that $ m_j\le ( \delta    t u ^2  ) ^{-1}$, and for this,
we can use the previous proof, with one additional fact.
Let $ \pi $ be a
finite partition of a probability space, and suppose that $ \pi '$ refines
$ \pi $.  Then, for any    random variable $ Z$ we have
\begin{equation*}
\mathbb E ( \mathbb E (Z\mid \pi ) ^2 )
\le \mathbb E ( \mathbb E (Z\mid \pi' ) ^2 )
\end{equation*}
Here, we are using a standard notation for conditional expectation given $ \pi $.
This is a simple martingale fact.
Indeed, let $ Y=\mathbb E (Z\mid \pi' )-\mathbb E (Z\mid \pi )$, and observe
that $ \mathbb E \; Y\cdot \mathbb E (Z\mid \pi )=0$.


The random variable in question is $ U_j$. Set
\begin{equation*}
S_{j,n}=\mathbb E \;\mathbb E (U_j\mid \mathcal P_n)^2 \,,\qquad 1\le n\le m\eqdef m_1+m_2+m_3\,.
\end{equation*}
We have just seen that the $ S _{j,n}$ are increasing in $1\le n\le m$.  
They are obviously at most $ \delta _{U_j}$.
And by construction, and in particular using (\ref{e.Up}) and (\ref{e.3Up}),
at each time at which the corresponding conditional is satisfied, we
increase these numbers by a definite amount.
At each iterate of the loop, in
a portion of the whole space $ H$ that is at least probability $ t \delta _{U_j}$,
where this uniformity constant is at least $ u$.
From this, we see that
\begin{equation*}
\delta _{U_j} \ge S _{j,\sigma _{j }(\ell )}\ge S _{j,\sigma _{j }(\ell -1)}+t u ^2 \delta _{U_j}\,,
\qquad 1\le \ell \le m _{j}
\,,\, j=1,2,3\,.
\end{equation*}
Therefore, $ m_j\le (t u ^2) ^{-1}  $.  The proof is complete.
\end{proof}

\begin{proof}[Proof of Lemma~\ref{uniformity}.]

Let us assume that e.\thinspace g.~$ D'=D$.
Apply Lemma~\ref{l.UU} to the sets $ X'=U_1$,  $ Y'=U_2$, $ D=D'=U_3$,
 with
$ t=\tfrac c {16} {\delta ^{2} }  $ and $ u=\upsilon ''$.
Let $  \mathcal P$ be the partition of $ H\times H$ that this Lemma gives us.

Define two subsets $ \mathcal E_1$ and $ \mathcal E_2$ of $ \mathcal P$ by
\begin{equation*}
\mathcal E_1 \eqdef \{ V_1 \times V_2 \in \mathcal P \mid  \mathbb P
(X'\cap V_j \mid V_j) \le t\}\,.
\end{equation*}
We define $ \mathcal E_2$ similarly, with $ X' $ replaced by $ Y'$,
and $ \mathcal E_3$ with $ X'$ replaced by $ D$.
These are the `empty' portions of the partition which nearly avoid $ X'$
or $ Y'$ entirely.

Consider
\begin{equation*}
X_0 \eqdef \bigcup \{X'\cap V_1\mid V_1 \times V_2\in  \mathcal E_1\cup \mathcal N_1\}
\end{equation*}
and similarly define $ Y_0$.  Then by
construction, $ \mathbb P (X_0\cap X')\le 2t=\frac c {8} {\delta ^{2} } $.
Observe that uniformity in $ D$ then implies that
\begin{equation*}
\mathbb E _{x,y\in H} X_0 (x) D (x+y) Y' (y)=\upsilon '+
\mathbb P (X_0\mid H) \delta _{D} \delta _{Y'}\,.
\end{equation*}
That is, we can assume
\begin{equation*}
\mathbb P (X_0 \times Y'\cap X_0 \stackrel {\textup{diag}} \times  D\mid X' \times Y'
\cap X' \stackrel {\textup{diag}} \times D)\le \tfrac c {4} {\delta ^{2} } \,.
\end{equation*}
The same inequality holds with $ X_0$ replaced by $ X$ and $ Y$ by $ Y_0$.
The import of this is the inequality
\begin{equation*}
\mathbb P (A\mid X_1 \times Y_1 \cap X_1 \stackrel {\textup{diag}} \times  D)
\ge \delta +\tfrac c2  {\delta ^{2} }\,, 	\qquad
X_1 \eqdef X-X_0\,,\, Y_1 \eqdef Y-Y_0\,.
\end{equation*}

Let $ V_1 \times V_2\in \mathcal P-\mathcal E_1 -\mathcal E_2-\mathcal N_1
- \mathcal N_2$.  In particular, $ X\cap V_1$ is a uniform set, obeying
$ \norm X\cap V_1. \textup{Uni}. \le \upsilon ''$.  As a consequence, we have
\begin{equation*}
\mathbb E _{x\in V_1, y\in V_2} X (x ) V_1 (x)
D (x+y) Y (y) V_2 (y)= \nu +\mathbb P (X\mid V_1) \mathbb P (D\mid V_1+V_2)
\mathbb P (Y\mid V_2)\,.
\end{equation*}
Here, $ \nu $ is a function of $ \upsilon '' $ that need not concern us.
Set
\begin{equation*}
D_0 \eqdef \{ D\cap V_1 +V_2 \mid V_1 \times V_2\in \mathcal E_3\cup \mathcal N_3\}.
\end{equation*}
It follows that
\begin{equation*}
\mathbb P (X_1 \times Y_1 \cap X_1 \stackrel {\textup{diag}} \times D_0
\mid X_1 \times Y_1 \cap X_1 \stackrel {\textup{diag}} \times D)
\le \tfrac c {4} {\delta ^{2 } } .
\end{equation*}

Then, for $ D_1=D-D_0$, we have
\begin{equation*}
\mathbb P (A\mid X_1 \times Y_1 \cap X_1 \stackrel {\textup{diag}} \times D_1)
\ge \delta +\tfrac c4 {\delta ^{2 } } .
\end{equation*}
From this, we see that we must have some $ H' \times H''\in \mathcal U$
for which
\begin{equation*}
\mathbb P (A\mid H' \times H''
\cap X_1 \times Y_1 \cap X_1 \stackrel {\textup{diag}} \times D_1)
\ge \delta +\frac c4  {\delta ^{2} } \,.
\end{equation*}
We take $ X''=X_1\cap V_1$, $ Y''=Y_1\cap V_2$, $ D''=D_1\cap
H'+H''$.  These are all uniform subsets, satisfying (\ref{conc1});
the second conclusion (\ref{conc2}) is the inequality above; the
lower bound on the dimension of $ H'$ and $ H''$ follows from Lemma~\ref{l.UU};
and the final conclusion follows from the fact that the element of the
partition that we chose, $ H' \times H''$ is not in the collection
$ \mathcal E_1\cup \mathcal E_2\cup \mathcal E_3$.

\end{proof}


\end{document}